\documentclass[a4paper]{article}
\usepackage{amsmath,amssymb,fancybox,longtable,graphics,amscd,multirow}
\usepackage[dvipdfm]{graphicx}
\usepackage{caption}
\newtheorem{thm}{\textsc{Theorem}}[section]
\newtheorem{prop}{\textsc{Proposition}}[section]
\newtheorem{lem}{\textsc{Lemma}}[section]
\newtheorem{remark}{\textsc{Remark}}
\newtheorem{cor}{\textsc{Corollary}}[section]
\newtheorem{defn}{\textsc{Definition}}[section]

\def\QED{$\Box$}
\def\END{$\blacksquare$}

\def\Pic{\mathop{\mathrm{Pic}}\nolimits}
\def\discr{\mathop{\mathrm{discr}}\nolimits}
\def\rk{\mathop{\mathrm{rk}}\nolimits}

\begin{document}
\title{A mirror duality for families of $K3$ surfaces associated to bimodular singularities}
\author{Makiko Mase}
\date{\small \begin{tabular}{l} Key Words: $K3$ surfaces, toric varieties, Picard lattices \end{tabular}}
\maketitle
\begin{abstract}
Ebeling and Ploog \cite{EbelingPloog} studied a duality of bimodular singularities which is part of the Berglund--H$\ddot{\textnormal{u}}$bsch mirror symmetry. 
Mase and Ueda \cite{MU} showed that this duality leads to a polytope mirror symmetry of families of $K3$ surfaces. 
We discuss in this article how this symmetry extends to a symmetry between lattices. 
\end{abstract}
\section{Introduction} \label{Introduction}
Bimodular singularities classified by Arnold \cite{Arnold75} have a duality studied by Ebeling and Ploog \cite{EbelingPloog} analogous to Arnold's strange duality for unimodular singularities. 
Namely, a pair $((B,f),\, (B',f'))$ of singularities $B,\, B'$ in $\mathbb{C}^3$ together with appropriate defining polynomials $f,\, f'$ is dual if the matrices $A_f,\, A_{f'}$ of exponents of $f$ and $f'$ are transpose to each other. 
Moreover, in some cases, such polynomials are compactified as anticanonical members of   $3$-dimensional weighted projective spaces whose general members are Gorenstein $K3$. 
The strange duality for unimodular singularities is related with the polytope mirror symmetry for families of $K3$ surfaces that are obtained by compactifying the singularities by Kobayashi \cite{Kobayashi} in a certain sense. 
In the study of bimodular singularities, Mase and Ueda \cite{MU} extend the duality by Ebeling and Ploog to a polytope mirror symmetry of families of $K3$ surfaces. 
More precisely, the following statement is shown : \\

\noindent
{\bf Theorem}\, \cite{MU} 
{\it Let $((B,\, f),\, (B',\, f'))$ be a dual pair in the sense of \cite{EbelingPloog} of  singularities $B$ and $B'$ together with their defining polynomials $f$ and $f'$  that are respectively compactified into polynomials $F$ and $F'$ as in \cite{EbelingPloog}. 
Then, there exists a reflexive polytope $\Delta$ such that $\Delta_F\subset \Delta$ and $\Delta_{F'}\subset \Delta^*$. 
Here, $\Delta_F$ and $\Delta_{F'}$ are respectively the Newton polytopes of $F$ and of $F'$, and $\Delta^*$ is the polar dual to $\Delta$. }
\END \\

In this article, we consider whether or not it is possible to extend the duality obtained in \cite{MU} further to the lattice mirror symmetry of families of $K3$ surfaces. 
More precisely, our problem is stated as follows: \\

\noindent
\begin{problem}\, 
Let $\Delta$ be a reflexive polytope as in \cite{MU}. 
Does there exist general members $S\in\mathcal{F}_\Delta$ and $S'\in\mathcal{F}_{\Delta^*}$ such that an isometry $T(\widetilde{S}) \simeq \Pic(\widetilde{S'})\oplus U$ holds ? 
Here, $\mathcal{F}_\Delta$ and $\mathcal{F}_{\Delta^*}$ are families of $K3$ surfaces associated to the polytopes $\Delta$ and $\Delta^*$, $\widetilde{S}$ denotes the minimal model of $S$, $\Pic(\widetilde{S})$ and $T(\widetilde{S})$ are respectively the Picard and transcendental lattice of $\widetilde{S}$. 
\end{problem} \\

The problem is answered in Theorem \ref{MainResult} together with an explicit description in Proposition \ref{ExplicitLattices} of the Picard lattices $\Pic(\Delta)$ and $\Pic(\Delta^*)$ defined in section \ref{MainThm}, with ranks $\rho(\Delta)$ and $\rho(\Delta^*)$, of the minimal model of appropriate general members in the families. 
The main result of this article is summarized here. 
In the sequel, the names of singularities follow Arnold \cite{Arnold75}, and singularities in a same row of Table \ref{NicePairsIntro} are dual to each other in the sense of \cite{EbelingPloog}. \\

\noindent
{\bf Proposition} \ref{ExplicitLattices} and {\bf Theorem} \ref{MainResult} \, 
{\it Let $\Delta$ be the reflexive polytope obtained in \cite{MU}. 
For the following transpose-dual pairs, the polar duality extends to a lattice mirror symmetry  between the families $\mathcal{F}_\Delta$ and $\mathcal{F}_{\Delta^*}$, where the Picard lattices are given in Table \ref{NicePairsIntro}. 
Here we use the notation $C^6_8:=\left( \begin{array}{cc} -4 & 1 \\ 1 & -2 \end{array}\right)$. }
\[
\begin{array}{cccccc} 
\textnormal{Singularity} & \Pic(\Delta) & \rho(\Delta) & \rho(\Delta^*) & \Pic(\Delta^*) & \textnormal{Singularity} \\ 
\hline \hline
Q_{12} & U\oplus E_6\oplus E_8 & 16 & 4 & U\oplus A_2 & E_{18}\\
\hline
Z_{1,0} & U\oplus E_7\oplus E_8 & 17 & 3 & U\oplus A_1 & E_{19}\\
\hline
E_{20} & U\oplus E_8^{\oplus 2} & 18 & 2 & U & E_{20}  \\
\hline
Q_{2,0} & U\oplus A_6\oplus E_8 & 16 & 4 & U\oplus C^6_8 & Z_{17}\\
\hline
E_{25} & U\oplus E_7\oplus E_8 & 17 & 3 & U\oplus A_1 & Z_{19}\\
\hline
Q_{18} & U\oplus E_6\oplus E_8 & 16 & 4 & U\oplus A_2 & E_{30}\\
\hline
\end{array}
\]
\begingroup
\captionof{table}{{\it Picard lattices for lattice mirror symmetric pairs}\END}\label{NicePairsIntro}
\endgroup
\noindent
\\

Section \ref{Preliminary} is to define the polytope- and lattice- mirror theories in subsection \ref{Mirrors}, and to define the transpose duality following \cite{EbelingPloog} in subsection \ref{EPDual}, based on a brief introduction to lattice theory in subsection \ref{Lattice}, and of toric geometry in subsection \ref{toricHypersurface}, where several formulas and results are stated without proof. 

The main theorem of this article is stated in section \ref{MainThm} following auxiliary results. 
The facts introduced in the previous section are used in their proof.

Denote by $\Delta_B$ the reflexive polytope obtained in \cite{MU} for a singularity $B$. 
As is seen in Table \ref{NicePairsIntro}, there are isometric Picard lattices  $\Pic(\Delta_{Q_{12}})\simeq \Pic(\Delta_{Q_{18}})$, and $\Pic(\Delta_{Z_{1,0}})\simeq \Pic(\Delta_{E_{25}})$. 
We consider and affirmatively answer in Proposition \ref{BirFamilies} the following question as an application in section \ref{Application}: \\

\noindent
\begin{problem}\, 
Are the families $\mathcal{F}_{\Delta_{Q_{12}}}$ (resp.\,$\mathcal{F}_{\Delta_{Z_{1,0}}}$) and $\mathcal{F}_{\Delta_{Q_{18}}}$ (resp.\,$\mathcal{F}_{\Delta_{E_{25}}}$) essentially the same in the sense that general members in these families are birationally equivalent ? 
\end{problem}  \\

\noindent
\begin{ackn}\\
{\rm {
The author thanks to the referee for his helpful comments particularly about the proof of the key Proposition \ref{ExplicitLattices} in the original manuscript. 
}}
\end{ackn}
\section{Preliminary} \label{Preliminary}
We start with having a consensus as to {\it Gorenstein $K3$ } and {\it $K3$ surfaces}. 
\begin{defn}
\rm
 A compact complex connected $2$-dimensional algebraic variety $S$ with at most $ADE$ singularities is called {\it Gorenstein $K3$} if (i)\, $K_S \sim 0$; and  (ii)\, $H^1(S,\, \mathcal{O}_S)=0$. 
If a Gorenstein $K3$ surface $S$ is nonsingular, $S$ is simply called a {\it $K3$ surface}. 
\END \end{defn}
\subsection{Brief lattice theory} \label{Lattice}
A {\it lattice} is a non-degenerate finitely-generated free $\mathbb{Z}$-module with a symmetric bilinear form called an intersection pairing. 
The {\it{discriminant group}} of a lattice $L$ is defined by $A_L:=L^*\slash L$, which is  finitely-generated and abelian, where $L^*:=Hom (L,\, \mathbb{Z})$ is dual to $L$. 
It is known that the order $|A_L|$ of the discriminant group is equal to the determinant of any intersection matrix of $L$.  
Let us recall a standard lattice theory by Nikulin \cite{Nikulin80}: 
\begin{cor}[Corollary 1.13.5-(1) \cite{Nikulin80}] \label{IncludeU}
If an even lattice $L$ of signature $(t_+,\, t_-)$ satisfies $t_+\geq 1,\, t_-\geq 1$, and $t_+ +  t_-\geq 3+{\rm length}\,  A_L$, then, there exists a lattice $T$ such that $L\simeq U\oplus T$, where $U$ is the hyperbolic lattice of rank $2$. \END
\end{cor}
In particular, if an even lattice $L$ is of $\rk L>12, t_+\geq 1,$ and $t_-\geq 1$, then, there exists a lattice $T$ such that $L\simeq U\oplus T$.  

Suppose $L$ is a sublattice of a lattice $L'$ with inclusion  $\iota: L\hookrightarrow L'$. 
Denote by $L^\perp_{L'}$ the orthogonal complement of $L$ in $L'$. 
The embedding $\iota$  is called {\it primitive}, and $L$ called a {\it primitive sublattice} of $L'$ if the finite abelian group $L'\, \slash \, \iota(L)$ is torsion-free. 
In other words, if there is no {\it overlattice} that is an intermediate lattice between $L$ and $L'$ of rank equal to the rank of $L$. 
Note that if a direct sum $L_1\oplus L_2$ of lattices is a sublattice of $L'$, then $(L_1\oplus L_2)^\perp_{L'} \simeq (L_1)^\perp_{L'}\oplus (L_2)^\perp_{L'}$. 

For a $K3$ surface $S$, it is known that $H^2(S,\, \mathbb{Z})$ with the intersection pairing is isometric to the {\it $K3$ lattice} $\Lambda_{K3}$ that is even unimodular of rank $22$ and signature $(3, 19)$, being $U^{\oplus 3}\oplus E_8^{\oplus 2}$, where $E_8$ is the negative-definite even unimodular lattice of rank $8$. 
There is a standard exact sequence $0\to H^1(S,\,\mathcal{O}_S^*)\overset{c_1}{\to}H^2(S,\,\mathbb{Z})\to 0$ of cohomologies so $H^1(S,\,\mathcal{O}_S^*)$ is inherited a lattice structure from $H^2(S,\,\mathbb{Z})$. 
Define the {\it Picard lattice} $\Pic(S)$ of a $K3$ surface $S$ as the group $c_1(H^{1,1}(S))\cap H^2(S,\,\mathbb{Z})$ with the lattice structure. 
The rank of $\Pic(S)$ is called the {\it Picard number}, denoted by $\rho(S)$. 
The Picard lattice is hyperbolic since a $K3$ surface is complex and algebraic, and is known to be a primitive sublattice of $\Lambda_{K3}$ under a {\it marking} $H^2(S,\, \mathbb{Z})\overset{\sim}{\to} \Lambda_{K3}$. 

If an even hyperbolic lattice $L$ of rank $\rk L \leq 20$ has $|A_L|$ being square-free, $L$ is a primitive sublattice of $\Lambda_{K3}$. 
Indeed, if so, for any lattice $L\subset L''\subset \Lambda_{K3}$ the general relation $|A_L| = [L'':L]^2 |A_{L''}|$ implies $[L'':L] = 1$ thus $L''\simeq L$. 
Hence there is no overlattice of $L$. 
Moreover, by surjectivity of the period mapping \cite{BHPV}, there exists a $K3$ surface $S$ such that $\Pic(S)\simeq L$. 
Let $M\subset \Lambda_{K3}$ be a hyperbolic sublattice. 
A $K3$ surface is {\it $M$-polarised} \cite{DolgachevMirror} if there exists a marking $\phi$ such that all divisors in $\phi^{-1}(C_M^{pol})$ are ample, where $C_M^{pol}$ is the positive cone in $M_\mathbb{R}$ minus $\bigcup_{d\in\Delta_M} H_d$, $\Delta_M=\{d\in M\, | \, d.d=-2\}$, and $H_d=\{ x\in \mathbb{P}(\Lambda_{K3})\, |\, x.d=0 \textnormal{ for all } d\in\Delta_M\}$. 

Nishiyama \cite{Nishiyama96} gives the orthogonal complements of primitive sublattices of type $ADE$ of $E_8$ in possible cases. 
\begin{lem}[Lemma 4.3 \cite{Nishiyama96}]\label{NLem4.3}
There exist primitive embeddings of lattices of type $ADE$ into $E_8$ with orthogonal complements given as follows. 
All the notation follows Bourbaki except $C_8^6 := \left( \begin{array}{cc} -4 & 1 \\ 1 & -2 \end{array}\right)$. 
\[
\begin{array}{llll}
(A_1)^\perp_{E_8}\simeq E_7 & (A_2)^\perp_{E_8}\simeq E_6 & (A_3)^\perp_{E_8}\simeq D_5 & (A_4)^\perp_{E_8}\simeq A_4 \\  (A_5)^\perp_{E_8}\simeq A_1\oplus A_2 & (A_6)^\perp_{E_8}\simeq C_8^6 & (A_7)^\perp_{E_8}\simeq ({-}8) \\ (D_4)^\perp_{E_8}\simeq D_4 & (D_5)^\perp_{E_8}\simeq A_3 & (D_6)^\perp_{E_8}\simeq A_1^{\oplus 2} & (D_7)^\perp_{E_8}\simeq ({-}4)  \\  (E_6)^\perp_{E_8}\simeq A_2 & (E_7)^\perp_{E_8}\simeq A_1 & & \hspace{20mm}\textnormal{\END}
\end{array}
\]
\end{lem}

\subsection{Brief toric geometry} \label{toricHypersurface}
Here we summarize toric divisors and $\Delta$-regularity. 
Let $M$ be a rank-$3$ lattice with the standard basis $\{ e_1,\, e_2,\, e_3\}$, $N$ be its dual, and $(\, ,\, ):M\times N \to \mathbb{Z}$ be the natural pairing. 
From now on, a polytope means a $3$-dimensional convex hull of finitely-many points in $\mathbb{Z}^3$ embedded into $\mathbb{R}^3$, namely, {\it integral}, and the origin is the only lattice point in the interior of it. 

Let $\mathbb{P}_\Delta$ be the toric variety defined by a polytope $\Delta$ in $M\otimes_\mathbb{Z}\mathbb{R}$, to which one can associate a fan $\Sigma_\Delta$ whose one-dimensional cones, called {\it one-simplices}, are generated by primitive lattice vectors each of whose end-point is an intersection point of $N$ and an edge of its polar dual $\Delta^*$ defined by 
\[
\Delta^* := \left\{ y\in N\otimes \mathbb{R} \, | \, (x,\, y) \geq -1  \textnormal{ for all } x\in \Delta \right\}. 
\] 
Let $\widetilde{\mathbb{P}_\Delta}$ be the toric resolution of singularities in $\mathbb{P}_\Delta$. 
A {\it{toric divisor}} is a divisor admitting the torus action, identified with the closure of the torus-action orbit of a one-simplex. 
Let ${\rm Div}_\mathbb{T}(\widetilde{\mathbb{P}_\Delta})$ be the set of all toric divisors $D_i=\overline{orb(\mathbb{R}_{\geq 0} v_i)},\, i=1,\ldots , s$ on $\widetilde{\mathbb{P}_\Delta}$, where $v_i$ is a primitive lattice vector, then, ${-}K_{\widetilde{\mathbb{P}_\Delta}}=\sum_{i=1}^s \, D_i$. 
By a standard exact sequence and a commutative diagram \cite{Oda78}
\[
\begin{array}{ccccccccc}
0 & \to & M & \to & {\rm Div}_\mathbb{T}(\widetilde{\mathbb{P}_\Delta}) & \to & \Pic(\widetilde{\mathbb{P}_\Delta}) & \to & 0 \\
 &  & &  & \downarrow & & \downarrow & &  \\
 &  & M & \to & \bigoplus_{i=1}^s \mathbb{Z}D_i & \to & A_2(\widetilde{\mathbb{P}_\Delta}) & \to & 0 \\
\end{array}
\]
there is a system of linear equations among toric divisors 
\begin{equation} \label{linearRelation}
\sum_{i=1}^s (e_j,\, v_i) D_i = 0, \hspace{3mm} j=1,2,3. 
\end{equation}
Thus the solution set of the linear system is generated by $(s-3)$ elements corresponding to  divisors which generate the Picard group $\Pic(\widetilde{\mathbb{P}_\Delta})$ of $\widetilde{\mathbb{P}_\Delta}$. 

\begin{defn}
\rm 
A polytope is {\it reflexive} if its polar dual is also integral. \END
\end{defn}

The importance of that we consider reflexive polytopes is by the following: 
\begin{thm}c.f.  \cite{BatyrevMirror}
The following conditions are equivalent: 
\begin{itemize}
\item[(1)] A polytope $\Delta$ is reflexive. 
\item[(2)] The toric $3$-fold $\mathbb{P}_\Delta$ is Fano, in particular, general anticanonical members of $\mathbb{P}_\Delta$ are Gorenstein $K3$. \END
\end{itemize}
\end{thm}

For a reflexive polytope $\Delta$, denote by $\mathcal{F}_{\Delta}$ the family of hypersurfaces parametrised by the complete anticanonical linear system of $\mathbb{P}_\Delta$. 
Note that general members in $\mathcal{F}_\Delta$ are Gorenstein $K3$ due to the previous theorem so that they are birationally equivalent to $K3$ surfaces by the existence of crepant resolution. 
Thus we call the family $\mathcal{F}_\Delta$ a {\it family of $K3$ surfaces}. 

We recall from \S 3 of \cite{BatyrevMirror} the notion of $\Delta$-regularity. 
\begin{defn} \rm
Let $F$ be a Laurent polynomial defining a hypersurface $Z_F$, whose Newton polytope is a polytope $\Delta$. 
The hypersurface $Z_F$ is called {\it $\Delta$-regular} if for every face $\Gamma$ of $\Delta$, the corresponding affine stratum $Z_{F, \Gamma}$ of $Z_F$ is either empty or a smooth subvariety of codimension $1$ in the torus $\mathbb{T}_\Gamma$ that is contained in the affine variety associated to $\Gamma$. \END
\end{defn}

It is shown in \cite{BatyrevMirror} that $\Delta$-regularity is a general condition, and singularities of all $\Delta$-regular members are simultaneously resolved by a toric desingularization of $\mathbb{P}_\Delta$. 
From now on, suppose a polytope $\Delta$ is reflexive and $S$ is a $\Delta$-regular member whose minimal model $\widetilde{S}$ is obtained by a toric resolution. 
\begin{defn} \label{DefOfL0}\rm
For a restriction $r:\widetilde{\mathbb{P}_\Delta}\to\widetilde{S}$, let $r_*:H^{1,1}(\widetilde{\mathbb{P}_\Delta})\to H^{1,1}(\widetilde{S})$ be the induced  mapping. 
Define a lattice $L_D(\widetilde{S}):=r_*(H^{1,1}(\widetilde{\mathbb{P}_\Delta}))\cap H^2(\widetilde{S},\, \mathbb{Z})$ of the intersection of the image of $r_*$ and $H^2(\widetilde{S},\,\mathbb{Z})$, and its orthogonal complement $L_0(\widetilde{S}):=L_D(\widetilde{S})^\perp_{H^2(\widetilde{S},\, \mathbb{Z})}$ in $H^2(\widetilde{S},\, \mathbb{Z})$. \END
\end{defn}

It is known \cite{Kobayashi} that $\rho(\widetilde{S})$ and $\rk L_0(\widetilde{S})$ only depend on the number of lattice points in edges of $\Delta$ and $\Delta^*$. 
Thus we define the Picard number $\rho(\Delta):=\rho(\widetilde{S})$, and the rank  $\rk L_{0,\Delta}:= \rk L_0(\widetilde{S})$ associated to $\Delta$. 
More precisely, denote by $\Gamma^*$ in $\Delta^*$ the dual face to a face $\Gamma$ of $\Delta$, and $l^*(\Gamma)$ is the number of lattice points in the interior of $\Gamma$, and $\Delta^{[1]}$ the set of edges in $\Delta$. 
Let $s$ be the number of one-simplices of $\Sigma_\Delta$. 
Then
\begin{eqnarray}
\rk L_{0,\Delta} & = & \sum_{\Gamma\in\Delta^{[1]}}l^*(\Gamma)\, l^*(\Gamma^*) = \rk L_{0,\Delta^*} , \label{rkL0} \\
\rho(\Delta) & = & s-3+\rk L_{0,\Delta}, \label{PicardNo} \\
\rho(\Delta) + \rho(\Delta^*) & = & 20 + \rk L_{0,\Delta}. \label{Pic+Trs} 
\end{eqnarray} 
If $l^*(\Gamma^*)=n_\Gamma$ and $l^*(\Gamma)=m_\Gamma$ for an edge $\Gamma$ of $\Delta$, there is a singularity of type $A_{n_\Gamma+1}$ with multiplicity $m_\Gamma+1$ on an affine variety associated to $\Gamma$. 

As we will see later, we only need formulas when $\rk L_{0,\Delta}=0$ for the intersection numbers of the divisors $\{ r_*D_i\}$ in $H^2(\widetilde{S},\,\mathbb{Z})$ given as follows. 
\begin{eqnarray}
r_*D_i^2 = r_*D_i.r_*D_i & = & \begin{cases} 2l^*(v_i^*)-2 & \textnormal{if $v_i$ is a vertex of $\Delta^*$, }\\ -2 & \textnormal{otherwise. }\end{cases} \label{SelfIntersection}
\end{eqnarray}
If end-points of $v_i$ and $v_j$ are on the edge $\Gamma_{ij}^*$ of $\Delta^*$, then 
\begin{eqnarray}
r_*D_i.r_*D_j & = & \begin{cases} 1 &  \textnormal{if $v_i$ and $v_j$ are next to each other, } \\ l^*(\Gamma_{ij})+1 & \textnormal{if $l^*(\Gamma_{ij}^*)=0$, and $v_i,\, v_j$ are both vertices,}  \\ 0 & \textnormal{otherwise. }\end{cases} \label{Intersection}
\end{eqnarray}

The Picard lattices of the minimal models of any $\Delta$-regular members, which are generated by components of restricted toric divisors, are isometric. 
Define the {\it{Picard lattice $\Pic(\Delta)$ of the family $\mathcal{F}_\Delta$}} as the Picard lattice of the minimal model of a $\Delta$-regular member with rank $\rho(\Delta)$. 
The orthogonal complement $T(\Delta)=\Pic(\Delta)^\perp_{\Lambda_{K3}}$ is called the {\it{transcendental lattice of $\mathcal{F}_\Delta$}}. 

\subsection{Mirrors} \label{Mirrors}
We define the polytope- and lattice-mirror theories. 
\subsubsection{Polytope Mirror}
We focus on {\it polytopes} that ``represent'' the anticanonical members in a toric variety as is seen in subsection \ref{toricHypersurface}. 

\begin{defn} \rm
A pair $(\Delta_1,\, \Delta_2)$ of reflexive polytopes or a pair $(\mathcal{F}_{\Delta_1},\, \mathcal{F}_{\Delta_2})$ of families of $K3$ surfaces associated to $\Delta_1$ and $\Delta_2$ is called {\it polytope mirror symmetric} if an isometry $\Delta_1 \simeq \Delta_2^*$ holds. \END
\end{defn}
\subsubsection{Lattice Mirror}
For a $K3$ surface $S$, $\Pic(S)$ is the Picard lattice, and $T(S)=\Pic(S)^\perp_{\Lambda_{K3}}$ is the transcendental lattice. 
A mirror for family of $M$-polarised $K3$ surfaces is defined when $M$ is a sublattice of $\Lambda_{K3}$ in general \cite{DolgachevMirror}. 
Here, we deal with the most strict case, namely, mirror for $K3$ surfaces with Picard lattice as their polarisation. 
\begin{defn} \rm 
(1)\, A pair $(S,\, S')$ of $K3$ surfaces is called {\it lattice mirror symmetric} if an isometry $T(S) \simeq \Pic(S') \oplus U$ holds.  \\
(2)\, A pair $(\mathcal{F},\, \mathcal{F}')$ of families whose general members are Gorenstein $K3$ surfaces is {\it lattice mirror symmetric} if there exist general members $S\in \mathcal{F}$ and $S'\in \mathcal{F}'$ pair of whose minimal models is lattice mirror symmetric. \END
\end{defn}

Note that a lattice mirror pair $(S,\, S')$ of $K3$ surfaces satisfies, by definition, $\rk \Pic(S') + 2 = \rk T(S) = 22-\rk\Pic(S)$, thus 
\begin{equation}
\rho(S')+\rho(S) = 20. \label{RhoDol}
\end{equation}

\subsection{Bimodular singularities and the transpose duality} \label{EPDual}
Being classified by Arnold \cite{Arnold75} in 1970's, bimodular singularities have two specific classes: quadrilateral and exceptional. 
Quadrilateral bimodular singularities are 6 in number with exceptional divisor of type $I_0^*$, whilst exceptional are 14 in number with exceptional divisor of type $II^*,\, III^*$ or $IV^*$ in Kodaira's notation. 

A non-degenerate polynomial $f$ in three variables is called {\it invertible} if $f$ has three terms  $f=\sum_{j=1}^3 x^{a_{1j}}y^{a_{2j}}z^{a_{3j}}$ such that its {\it matrix} $A_f:=(a_{ij})_{1\leq i,j\leq 3}$ {\it of exponents} is invertible in $GL_3(\mathbb{Q})$.

\begin{defn}c.f. \rm \cite{EbelingPloog} 
Let $B=(0,(f=0))$ and $B'=(0,(f'=0))$ be germs of singularities in $\mathbb{C}^3$. 
A pair $(B, B')$ of singularities is called {\it transpose dual} if the following three conditions are satisfied. 
\begin{itemize}
\item[(1)] Defining polynomials $f$ and $f'$ are invertible. 
\item[(2)] Matrices $A_f$ and $A_{f'}$ of exponents of $f$ and $f'$ are transpose to each other. 
\item[(3)] $f$ $($resp. $f')$ is compactified to a four-term polynomial $F$ $($resp. $F')$ in $|{-}K_{\mathbb{P}(a)}|$ $($resp. $|{-}K_{\mathbb{P}(b)}|)$, where $\mathbb{P}(a)$ $($resp. $\mathbb{P}(b))$ is the $3$-dimensional weighted projective space whose general members are Gorenstein $K3$ with weight $a$ $($resp. $b)$ out of the list of $95$ weights classified by \cite{Yonemura}\cite{Iano-Fletcher}\cite{C3f}. \END
\end{itemize}
\end{defn}

Condition (1) and (2) is said that they are Berglund--H$\ddot{\textnormal{u}}$bsch mirror symmetric.  
Ebeling-Ploog \cite{EbelingPloog} show that there are 16 transpose-dual pairs among quadrilateral and exceptional bimodular singularities, exceptional unimodular singularities, and the singularities $E_{25},\, E_{30},\, X_{2,0}$, and $Z_{2,0}$. 

\section{The transpose dual and the lattice mirror} \label{MainThm}
For a transpose-dual pair $(B, B')$ of defining polynomial $f$ (resp. $f'$) being compactified to a polynomial $F$ (resp. $F'$), consider the Newton polytope $\Delta_F$ (resp. $\Delta_{F'}$) of $F$ (resp. $F'$) all of whose corresponding monomials are fixed by an automorphism action on $(F=0)$ (resp. $(F'=0)$). 
\begin{remark} \rm
A compactified member $F$ to $f$ does not always define a Gorenstein $K3$ surface because $\Delta_F$ may not be reflexive. 
\end{remark}

However, the Newton polytopes are extended to be reflexive and dual. 

\begin{thm}\cite{MU} \label{MU}
For each transpose-dual pair $(B, B')$, there exists a reflexive polytope $\Delta$ such that $\Delta_F\subset \Delta$ and $\Delta_{F'}\subset \Delta^*$. \END
\end{thm}

Computing Picard lattices is generally difficult, but it seems possible for $\Delta$-regular members by subsection \ref{toricHypersurface}. 
Let us reformulate our problem. \\

\noindent
\begin{problem}\, 
For a polytope $\Delta$ obtained in \cite{MU}, is a pair $(\widetilde{S},\, \widetilde{S'})$ of minimal models of $\Delta$-regular $S\in\mathcal{F}_\Delta$ and $\Delta^*$-regular $S'\in\mathcal{F}_{\Delta^*}$ lattice mirror symmetric ?
\end{problem} \\

First we study the rank $\rk L_{0,\Delta}$. 
\begin{lem} \label{RkOfL0}
The list of $\rk L_{0, \Delta}$ for the reflexive polytope $\Delta$ obtained in \cite{MU} is given in Table \ref{RankL0}. 
\[
\begin{array}{ccccc}  \label{RankL0Lemma}
\textnormal{transpose-dual pair} & \rk L_{0, \Delta} & & \textnormal{transpose-dual pair} & \rk L_{0, \Delta} \\ 
\hline
\hline
(Q_{12},\, E_{18}) & 0 & & (Z_{17},\, Q_{2,0}) & 2 \\
\hline
(Z_{1,0},\, E_{19}) & 0 & & (U_{1,0},\, U_{1,0}) & 2 \\
\hline
(E_{20},\, E_{20}) & 0 & & (U_{16},\, U_{16}) & 2 \\
\hline
(Q_{2,0},\, Z_{17}) & 0 & & (Q_{17},\, Z_{2,0}) & 2 \\
\hline
(E_{25},\, Z_{19}) & 0 & & (W_{1,0},\, W_{1,0}) & 3 \\
\hline
(Q_{18},\, E_{30}) & 0 & & (W_{17},\, S_{1,0}) & 5 \\
\hline
(Z_{1,0},\, Z_{1,0}) & 1 & & (W_{18},\, W_{18}) & 6 \\
\hline
(Z_{13},\, J_{3,0}) & 2 & &  (S_{17},\, X_{2,0}) & 6 \\
\hline 
\end{array}
\]
\captionof{table}{$\rk L_{0, \Delta}$}\label{RankL0}
\end{lem}
{\sc Proof. }\, 
The assertion follows from direct and case-by-case computation by formula (\ref{rkL0}) in subsection \ref{toricHypersurface}. \\
{\bf 1.}\, $(Q_{12},\, E_{18})$ \, 
The polytope $\Delta$ is given in subsection 4.7 of \cite{MU}. 
There is no contribution to $\rk L_{0,\Delta}$ since $l^*(\Gamma)$ or $l^*(\Gamma^*)$  is zero for any edge $\Gamma$. 
Thus, $\rk L_{0, \Delta} = 0$ by formula (\ref{rkL0}). \\
Similar for the cases $(Z_{1,0},\, E_{19}),\, (E_{20},\, E_{20}),\, (Q_{2,0},\, Z_{17}),\, (E_{25},\, Z_{19}) ,\, (Q_{18},\, E_{30})$. \\
{\bf 2.}\, $(Z_{13},\, J_{3,0})$ \, 
The polytope $\Delta$ is given in subsection 4.1 of \cite{MU}. 
The only contribution to $\rk L_{0, \Delta}$ is by an edge $\Gamma$ between vertices $(0,\, 0,\, 1)$ and $(-2,\, -6,\, -9)$, whose dual $\Gamma^*$ is between $(8,\, -1,\, -1)$ and $(-1,\, 2,\, -1)$ so $l^*(\Gamma)=1$ and $l^*(\Gamma^*)=2$. 
Thus, $\rk L_{0, \Delta} = 1\times 2 = 2$ by formula (\ref{rkL0}). \\
Similar for the cases $(Z_{1,0}, Z_{1,0}),\, (Z_{17}, Q_{2,0}),\, (U_{1,0}, U_{1,0}),\, (U_{16}, U_{16}),\, (Q_{17}, Z_{2,0}),$\\
$(W_{1,0}, W_{1,0}),\, (W_{17}, S_{1,0}),\, (W_{18}, W_{18}),\, (S_{17}, X_{2,0})$. 
\QED
\begin{cor}  \label{NotAdmitDolMir}
No $\Delta$- and $\Delta^*$-regular members for transpose-dual pairs $(Z_{13}, J_{3,0}),$\, $(Z_{1,0}, Z_{1,0}),$\, $(Z_{17}, Q_{2,0}),$\, $(U_{1,0}, U_{1,0}),$ \,
$(U_{16}, U_{16}),$\, $(Q_{17}, Z_{2,0}),$\, $(W_{1,0}, W_{1,0}),$\, $(W_{17}, S_{1,0}),$\, $(W_{18}, W_{18}),$\, $(S_{17}, X_{2,0})$ admit a lattice mirror symmetry. 
\end{cor}
{\sc Proof. }\, 
For each $\Delta$ associated to the presented pairs, by formula (\ref{Pic+Trs}),  
\[
\rho(\Delta) + \rho(\Delta^*) = 20 + \rk L_{0, \Delta} > 20
\]
since $\rk L_{0, \Delta} > 0$ by Lemma \ref{RankL0Lemma}. 
Thus, the equation (\ref{RhoDol}) does not hold. 
Therefore, $\Delta$- and $\Delta^*$-regular members do not admit a lattice mirror symmetry. 
\QED
\begin{cor} \label{SurjectiveRestriction}
The restriction mapping $r_* : \Pic(\widetilde{\mathbb{P}_\Delta}) \to \Pic(\widetilde{S})$, for $\Delta$-regular $S\in\mathcal{F}_\Delta$ is surjective for the transpose-dual pairs 
$(Q_{12}, E_{18}),$\, $(Z_{1,0}, E_{19}),$\, $(E_{20}, E_{20}),$\, $(Q_{2,0}, Z_{17}),$\, $(E_{25}, Z_{19}),$\, $(Q_{18}, E_{30}). $
\end{cor}
{\sc Proof.}\, 
By Lemma \ref{RankL0Lemma},  $\rk L_{0, \Delta} = 0$ for each case. 
By definition, $\rk L_{0, \Delta}$ is equal to the rank of the orthogonal complement of $r_*(H^{1,1}(\widetilde{\mathbb{P}_\Delta}))$ in $H^2(\widetilde{S},\, \mathbb{Z})$. 
Thus, that $\rk L_{0, \Delta} = 0$ means that $r_*$ is surjective. \QED \\

By Corollary \ref{NotAdmitDolMir}, we may only focus on the transpose-dual pairs appearing in Corollary \ref{SurjectiveRestriction}, whose statement means moreover that $\Pic(\Delta)$ is generated by restricted toric divisors generating $\Pic(\widetilde{\mathbb{P}_\Delta})$, and analogous to $\Delta^*$. \\

Let $A_L$ denote the discriminant group, $q_L$ the quadratic form, and $\discr L$ the discriminant of a lattice $L$. 
If $p=\discr L$ is prime, then $A_L\simeq \mathbb{Z}\slash p\mathbb{Z}$. 
Before stating our main results, note a fact in Proposition 1.6.1 in \cite{Nikulin80}. 
Suppose that lattices $S$ and $T$ are primitively embedded into the $K3$ lattice $\Lambda_{K3}$. 
If $A_S\simeq A_T$ and $q_S = -q_T$, then, it is determined that the orthogonal complement $S^\perp_{\Lambda_{K3}}$ in $\Lambda_{K3}$ is $T$. 
And $q_S=-q_T$ if and only if $\discr S = -\discr T$. 


\begin{prop} \label{ExplicitLattices}
The Picard lattice $\Pic(\Delta^*)$ for $\Delta$ in Corollary \ref{SurjectiveRestriction} is as in Table \ref{NicePairs}, where singularities in a row are transpose-dual. 
In each case, one gets 
\[
\discr \Pic(\Delta)=-\discr (U\oplus \Pic(\Delta^*)), \textnormal{ and } A_{\Pic(\Delta)}\simeq A_{U\oplus \Pic(\Delta^*)}. 
\]
Denote $C_8^6 := \left( \begin{array}{cc} -4 & 1 \\ 1 & -2 \end{array}\right)$. 
\[
\begin{array}{ccccc} 
\textnormal{Singularity} & \rho(\Delta^*) & \Pic(\Delta^*) & \textnormal{Singularity} \\ 
\hline
\hline
Q_{12} &  4 & U\oplus A_2 & E_{18}\\
\hline
Z_{1,0} & 3 & U\oplus A_1 & E_{19}\\
\hline
E_{20} & 2 & U & E_{20}  \\
\hline
Q_{2,0} & 4 & U\oplus C_8^6  & Z_{17}\\
\hline
E_{25} & 3 & U\oplus A_1 & Z_{19}\\
\hline
Q_{18} & 4 & U\oplus A_2 & E_{30}\\
\hline
\end{array}
\]
\begingroup
\captionof{table}{$\Pic(\Delta^*)$ for $\Delta$ in Corollary \ref{SurjectiveRestriction}}\label{NicePairs}
\endgroup
\end{prop}

\noindent
{\sc Proof.} \\ 
${\bf 1.}\, Q_{12}$ and $E_{18}$ \, 
The polytope $\Delta$ is given in subsection 4.7 of \cite{MU} and the associated toric $3$-fold has $19$ toric divisors $D_i$ corresponding to the one-simplices generated by vectors
\[
\begin{array}{lll}
m_1=(1,\,2,\,2) & m_2=(0,\,1,\,1) & m_3=(-8,\,-11,\,-9) \\
m_4=(1,\,-2,\,0) & m_5=(1,\,1,\,0) & m_6=(-4,\,-5,\,-4) \\
m_7=(-7,\,-10,\,-8) & m_8=(-6,\,-9,\,-7) & m_9=(-5,\,-8,\,-6) \\
m_{10}=(-4,\,-7,\,-5) & m_{11}=(-3,\,-6,\,-4) & m_{12}=(-2,\,-5,\,-3) \\
m_{13}=(-1,\,-4,\,-2) & m_{14}=(0,\,-3,\,-1) & m_{15}=(1,\,0,\,1) \\
m_{16}=(-5,\,-7,\,-6) & m_{17}=(-2,\,-3,\,-3) & m_{18}=(1,\,-1,\,0) \\
m_{19}=(1,\,0,\,0)
\end{array}
\]
and by solving the linear system (\ref{linearRelation}) : $\sum_{i=1}^{19}(e_j,\, m_i)D_i = 0 \quad (j=1,2,3)$, we get linear relations among toric divisors
\begin{eqnarray*}
D_1 & \sim & -9 D_4+3 D_5-D_7-2 D_8-3 D_9 -4 D_{10}-5 D_{11}-6 D_{12}-7 D_{13}-8 D_{14}\\
 & & \hspace{53mm} -5 D_{15}+D_{16}+2 D_{17}-5 D_{18}-D_{19}, \\
D_2 & \sim & D_3+10 D_4-2 D_5+2 D_7+3 D_8+4 D_9+ 5 D_{10}+6 D_{11}+7 D_{12}+8 D_{13} \\
 & & \hspace{50mm} +9 D_{14}+5 D_{15}-D_{17}+6D_{18}+2 D_{19}, \\
D_6 & \sim & -2 D_3-2 D_4+D_5-2 D_7-2 D_8-2 D_9 -2 D_{10}-2 D_{11}-2 D_{12}-2 D_{13} \\
 & & \hspace{65mm} -2 D_{14}-D_{15}-D_{16}-D_{18}. 
\end{eqnarray*}
So the set $\{D_i \, | \, i\not= 1,2,6 \}$ of toric divisors is linearly independent. 
Let $L$ be the lattice generated by the set $\{r_*D_i \, | \, i\not= 1,2,6 \}$ of their restrictions to a $\Delta$-regular member. 
We shall check that $L$ is primitively embedded into the $K3$ lattice to show that $L$ is indeed the Picard lattice of the family $\mathcal{F}_\Delta$. 
By computer calculation with formulas (\ref{SelfIntersection}) and (\ref{Intersection}), the determinant of an intersection matrix of $L$ is $-3$ since this matrix is given by 
\[
\left(
\begin{array}{cccccccccccccccc}
 -2 & 0 & 0 & 1 & 0 & 0 & 0 & 0 & 0 & 0 & 0 & 0 & 1 & 0 & 0 & 0 \\
 0 & -2 & 0 & 0 & 0 & 0 & 0 & 0 & 0 & 0 & 1 & 1 & 0 & 0 & 1 & 0 \\
 0 & 0 & 0 & 0 & 0 & 0 & 0 & 0 & 0 & 0 & 0 & 0 & 0 & 1 & 0 & 1 \\
 1 & 0 & 0 & -2 & 1 & 0 & 0 & 0 & 0 & 0 & 0 & 0 & 0 & 0 & 0 & 0 \\
 0 & 0 & 0 & 1 & -2 & 1 & 0 & 0 & 0 & 0 & 0 & 0 & 0 & 0 & 0 & 0 \\
 0 & 0 & 0 & 0 & 1 & -2 & 1 & 0 & 0 & 0 & 0 & 0 & 0 & 0 & 0 & 0 \\
 0 & 0 & 0 & 0 & 0 & 1 & -2 & 1 & 0 & 0 & 0 & 0 & 0 & 0 & 0 & 0 \\
 0 & 0 & 0 & 0 & 0 & 0 & 1 & -2 & 1 & 0 & 0 & 0 & 0 & 0 & 0 & 0 \\
 0 & 0 & 0 & 0 & 0 & 0 & 0 & 1 & -2 & 1 & 0 & 0 & 0 & 0 & 0 & 0 \\
 0 & 0 & 0 & 0 & 0 & 0 & 0 & 0 & 1 & -2 & 1 & 0 & 0 & 0 & 0 & 0 \\
 0 & 1 & 0 & 0 & 0 & 0 & 0 & 0 & 0 & 1 & -2 & 0 & 0 & 0 & 0 & 0 \\
 0 & 1 & 0 & 0 & 0 & 0 & 0 & 0 & 0 & 0 & 0 & -2 & 0 & 0 & 0 & 0 \\
 1 & 0 & 0 & 0 & 0 & 0 & 0 & 0 & 0 & 0 & 0 & 0 & -2 & 1 & 0 & 0 \\
 0 & 0 & 1 & 0 & 0 & 0 & 0 & 0 & 0 & 0 & 0 & 0 & 1 & -2 & 0 & 0 \\
 0 & 1 & 0 & 0 & 0 & 0 & 0 & 0 & 0 & 0 & 0 & 0 & 0 & 0 & -2 & 1 \\
 0 & 0 & 1 & 0 & 0 & 0 & 0 & 0 & 0 & 0 & 0 & 0 & 0 & 0 & 1 & -2
\end{array}
\right), 
\]
Since the discriminant of $L$ is $-3$ that is square-free, there exists no overlattice of $L$ ;  indeed, if  $H\subset \Lambda_{K3}$ were an overlattice of $L$, then, by the standard relation $-3=[H:L]^2 \discr H$, we get $[H:L]=1$ and  $\discr H=-3$ so that $L\simeq H$. 
Hence, $L$ is primitively embedded into the $K3$ lattice. 
By construction, $L$ is indeed the Picard lattice of the family $\mathcal{F}_\Delta$. 

The dual polytope $\Delta^*$ associates a toric $3$-fold with $7$ toric divisors $D'_i$ corresponding to the one-simplices generated by vectors 
\[
\begin{array}{llll}
v_1=(1,\,1,\,-2) & v_2=(1,\,-2,\,1) & v_3=(2,\,-3,\,2) & v_4=(-1,\,0,\,1) \\
v_5=(-1,\,0,\,0) & v_6=(1,\,0,\,-1) & v_7=(1,\,-1,\,0) 
\end{array}
\]
and by solving the linear system (\ref{linearRelation}) : $\sum_{i=1}^7(e_j,\, m_i)D'_i = 0\quad (j=1,2,3)$, we get linear relations among toric divisors
\[
D'_1 \sim  2 D'_2+3 D'_3+D'_7, \quad 
D'_4 \sim  3 D'_2+4 D'_3+D'_6+2 D'_7, \quad
D'_5 \sim  D'_3 . 
\]
Thus the set $\{D'_i \, | \, i\not= 1,4,5 \}$ of toric divisors is linearly independent. 
Let $L'$ be the lattice generated by the set $\{r_*D'_i \, | \, i\not= 1,4,5 \}$ of their restrictions to a $\Delta^*$-regular member. 
By formulas (\ref{SelfIntersection}) and (\ref{Intersection}), an intersection matrix associated to $L'$ is given by 
\[
\left(
\begin{array}{cccc}
 -2 & 1 & 0 & 1 \\
 1 & 0 & 0 & 0 \\
 0 & 0 & -2 & 1 \\
 1 & 0 & 1 & -2
\end{array}
\right)
\]
that is equivalent to $U\oplus A_2$ by re-taking the generators as
\[
\{r_*\, D'_2+r_*\, D'_3,\, r_*\, D'_3,\, r_*\, D'_6, \, r_*\, D'_3 + r_*\, D'_7\}. 
\]
Since the discriminant of $L'$ is $-3$ that is square-free, there exists no overlattice of $L'$ ;  indeed, if $H'\subset \Lambda_{K3}$ were an overlattice of $L'$, then, by the standard relation $-3=[H':L']^2 \discr H'$, we get $[H':L']=1$ and  $\discr H'=-3$ so that $H'\simeq L'$. 
Hence, $L'$ is primitively embedded into the $K3$ lattice. 
By construction, $L'$ is indeed the Picard lattice of the family $\mathcal{F}_{\Delta^*}$. 
Therefore $\Pic(\Delta^*)=L'=(\mathbb{Z}^4,\, U\oplus A_2)$. 

Similarly, the lattice $U\oplus \Pic(\Delta^*)$ is also primitively embedded into the $K3$ lattice since it is of signature $(2,\, 4)$ and is of discriminant $3$. 
Besides, $\discr(U\oplus\Pic(\Delta^*))=3=-\discr \Pic(\Delta)$, and moreover,  $A_{\Pic(\Delta)}\simeq A_{U\oplus\Pic(\Delta^*)}\simeq\mathbb{Z}\slash 3\mathbb{Z}$. \\

\noindent
${\bf 2. }\, Z_{1,0}$ and $E_{19}$ \, 
The polytope $\Delta$ is given in subsection 4.8 of \cite{MU} and the associated toric $3$-fold has $20$ toric divisors $D_i$ corresponding to 
\[
\begin{array}{lll}
m_1=(1,\,-1,\,2) & m_2=(0,\,-1,\,1) & m_3=(0,\,0,\,1) \\
m_4=(4,\,2,\,-1) & m_5=(-6,\,2,\,-11) & m_6=(2,\, 0 ,\, 1) \\
m_7=(3,\, 1,\, 0) & m_8=(-2,\, 0,\, -3) & m_9=(-4,\, 1,\, -7) \\
m_{10}=(2,\,1,\,0) & m_{11}=(-3,\,1,\,-5) & m_{12}=(3,\,2,\,-2) \\
m_{13}=(2,\,2,\,-3) & m_{14}=(1,\,2,\,-4) & m_{15}=(0,\,2,\,-5) \\
m_{16}=(-1,\,2,\,-6) & m_{17}=(-2,\,2,\,-7) & m_{18}=(-3,\,2,\,-8) \\
m_{19}=(-4,\,2,\,-9) & m_{20}=(-5,\, 2,\, -10)
\end{array}
\]
and $\{D_i\, | \, i\not=1,2,3 \}$ is linearly independent by system (\ref{linearRelation}). 
The lattice $L:=\langle r_* D_i\, | \, i\not=1,2,3 \rangle_\mathbb{Z}$ has $\rk L =17$, and $\discr L=2$ by an explicit calculation of its intersection matrix using formulas (\ref{SelfIntersection}) and (\ref{Intersection}).
Since the discriminant of $L$ is square-free, $L$ is primitively embedded into the $K3$ lattice, and thus $L=\Pic(\Delta)$. 

The dual polytope $\Delta^*$ associates a toric $3$-fold with $6$ toric divisors $D'_i$ corresponding to 
\[
\begin{array}{lll}
v_1=(1,\,-3,\,-1) & v_2=(2,\,0,\,-1) & v_3=(1,\,0,\,-1) \\
v_4=(0,\,-1,\,-1) & v_5=(-1,\,2,\,1) & v_6=(0,\,1,\,0) 
\end{array}
\]
and $\{D'_i\, |\, i\not=1,2,5\}$ is linearly independent by system (\ref{linearRelation}). 
By re-taking the generators as $\{ r_*\, D'_3+r_*\, D'_4,\, r_*\, D'_4,\, r_*\, D'_6-r_*\, D'_4\}$, the lattice $L':=\langle r_*D'_i \, | \, i\not= 1,2,5 \rangle_\mathbb{Z}$ has an intersection matrix $U\oplus A_1$. 
Since $\rk L' = 3$, and $\discr L' = 2$, which is square-free, $L'$ is primitively embedded into the $K3$ lattice, and thus $\Pic(\Delta^*)=L'=(\mathbb{Z}^3,\, U\oplus A_1)$. 

Similarly $U\oplus \Pic(\Delta^*)$ is primitively embedded into the $K3$ lattice.
Besides, $\discr (U\oplus \Pic(\Delta^*))=-2=-\discr \Pic(\Delta)$, and moreover,  $A_{\Pic(\Delta)}\simeq A_{U\oplus\Pic(\Delta^*)}\simeq\mathbb{Z}\slash 2\mathbb{Z}$. \\

\noindent
${\bf 3. }\, E_{20}$ and $E_{20}$ \, 
The polytope $\Delta$ is given in subsection 4.9 of \cite{MU} and the associated toric $3$-fold has $21$ toric divisors $D_i$ corresponding to 
\[
\begin{array}{lll}
m_1=(-1,\,-1,\,2) & m_2=(-1,\,-1,\,-1) & m_3=(-1,\,11,\,2) \\
m_4=(1,\,-1,\,0) & m_5=(-1,\,-1,\,1) & m_6=(-1,\, -1 ,\, 0) \\
m_7=(-1,\, 3,\, 0) & m_8=(-1,\, 7,\, 1) & m_9=(0,\, -1,\, 1) \\
m_{10}=(0,\,5,\,1) & m_{11}=(-1,\, 0,\,2) & m_{12}=(-1,\, 1,\,2) \\
m_{13}=(-1,\, 2,\,2) & m_{14}=(-1,\, 3,\,2) & m_{15}=(-1,\, 4,\,2) \\
m_{16}=(-1,\, 5,\,2) & m_{17}=(-1,\, 6,\,2) & m_{18}=(-1,\, 7,\,2) \\
m_{19}=(-1,\, 8,\,2) & m_{20}=(-1,\, 9,\, 2) & m_{21}=(-1,\, 10,\, 2)\end{array}
\]
and $\{D_i\, | \, i\not=10,13,14 \}$ is linearly independent by system (\ref{linearRelation}). 
The lattice $L:=\langle r_* D_i\, | \, i\not=10,13,14 \rangle_\mathbb{Z}$ has $\rk L = 18$, and $\discr L=-1$ by an explicit calculation of its intersection matrix using formulas (\ref{SelfIntersection}) and (\ref{Intersection}).
As it being unimodular, $L$ is primitively embedded into the $K3$ lattice, and thus $L=\Pic(\Delta)$. 

The dual polytope $\Delta^*$ associates a toric $3$-fold with $5$ toric divisors $D'_i$ corresponding to 
\[
\begin{array}{lll}
v_1=(0,\,1,\,0) & v_2=(1,\,0,\,0) & v_3=(-1,\,0,\,-1) \\
v_4=(-2,\,-1,\,4) & v_5=(-1,\,0,\,2)  
\end{array}
\]
and $\{D'_i\, |\, i\not=1,2,3\}$ is linearly independent by system (\ref{linearRelation}). 
By an intersection matrix computed by formulas (\ref{SelfIntersection}) and (\ref{Intersection}), the lattice $L':=\langle r_*D'_i \, | \, i\not= 1,2,3 \rangle_\mathbb{Z}$ has $\discr L' = -1,\, \rk L'=2$, and $L'$ is even. 
By the classification of even unimodular lattices, $L'$ is isometric to $U$, which is primitively embedded into the $K3$ lattice, and thus $\Pic(\Delta^*)=L'=(\mathbb{Z}^2,\, U)$. 

Similarly, $U\oplus\Pic(\Delta^*)$ is primitively embedded into the $K3$ lattice. 
Besides, $\discr (U\oplus\Pic(\Delta^*))=1=-\discr\Pic(\Delta)$, and moreover,  $A_{\Pic(\Delta)}\simeq A_{U\oplus\Pic(\Delta^*)}\simeq \{ 0\}$. \\

\noindent
${\bf 4. }\, Q_{2,0}$ and $Z_{17}$ \, 
The polytope $\Delta$ is given in subsection 4.10 of \cite{MU} and the associated toric $3$-fold has $19$ toric divisors $D_i$ corresponding to 
\[
\begin{array}{lll}
m_1=(0,\,1,\,1) & m_2=(1,\,2,\,2) & m_3=(1,\,1,\,2) \\
m_4=(0,\,-1,\,0) & m_5=(-6,\,-7,\,-9) & m_6=(1,\, 0 ,\, -2) \\
m_7=(1,\, 0,\, 1) & m_8=(1,\, 1,\, 0) & m_9=(-2,\, -3,\, -3) \\
m_{10}=(-4,\, -5,\, -6) & m_{11}=(1,\, 0,\, -1) & m_{12}=(1,\, 0,\, 0) \\
m_{13}=(-5,\, -6,\,-8) & m_{14}=(-4,\, -5,\,-7) & m_{15}=(-3,\, -4,\,-6) \\
m_{16}=(-2,\, -3,\,-5) & m_{17}=(-1,\, -2,\,-4) & m_{18}=(0,\, -1,\,-3) \\
m_{19}=(-3,\, -3,\,-4) 
\end{array}
\]
and $\{D_i\, | \, i\not=1,2,3 \}$ is linearly independent by system (\ref{linearRelation}). 
The lattice $L:=\langle r_* D_i\, | \, i\not=1,2,3 \rangle_\mathbb{Z}$ has $\rk L = 16$, and $\discr L=-7$ by an explicit calculation of its intersection matrix using formulas (\ref{SelfIntersection}) and (\ref{Intersection}). 
Since the discriminant of $L$ is square-free, $L$ is primitively embedded into the $K3$ lattice, and thus $L=\Pic(\Delta)$. 

The dual polytope $\Delta^*$ associates a toric $3$-fold with $7$ toric divisors $D'_i$ corresponding to 
\[
\begin{array}{llll}
v_1=(1,\,-2,\,1) & v_2=(-1,\,1,\,0) & v_3=(2,\,1,\,-2) & v_4=(1,\,0,\,-1) \\
v_5=(0,\,1,\,-1)  & v_6=(-1,\,0,\,0) & v_7=(1,\,-1,\,0)
\end{array}
\]
and $\{D'_i\, |\, i\not=1,2,6\}$ is linearly independent by system (\ref{linearRelation}).
By re-taking the generators as $\{r_*\, D'_3,\, r_*\, D'_3+r_*\,D'_4,\, 2r_*\, D'_3+r_*\, D'_4-r_*\, D'_5, \,r_*\, D'_7-r_*\, D'_3\}$, the lattice $L':=\langle r_*D'_i \, | \, i\not= 1,2,6 \rangle$ has an intersection matrix $U\oplus C_8^6$. 
Since $\discr L' = -7$ is square-free, $L'$ is primitively embedded into the $K3$ lattice, and thus $\Pic(\Delta^*)=L'=\left(\mathbb{Z}^4,\, U\oplus C_8^6\right)$. 

Similarly $U\oplus\Pic(\Delta^*)$ is primitively embedded into the $K3$ lattice. 
Besides, $\discr (U\oplus \Pic(\Delta^*)) = 7=-\discr \Pic(\Delta)$, and moreover, $A_{\Pic(\Delta)}\simeq A_{U\oplus\Pic(\Delta^*)}\simeq \mathbb{Z}\slash 7\mathbb{Z}$. \\

\noindent
${\bf 5. }\, E_{25}$ and $Z_{19}$ \, 
The polytope $\Delta$ is given in subsection 4.11 of \cite{MU} and the associated toric $3$-fold has $20$ toric divisors $D_i$ corresponding to 
\[
\begin{array}{lll}
m_1=(-1,\,2,\,0) & m_2=(-1,\,-1,\,9) & m_3=(-1,\,-1,\,-1) \\
m_4=(-1,\,2,\,-1) & m_5=(1,\,-1,\,-1) & m_6=(-1,\, 1 ,\, 3) \\
m_7=(-1,\, 0,\, 6) & m_8=(-1,\, -1,\, 8) & m_9=(-1,\, -1,\, 7) \\
m_{10}=(-1,\,-1,\, 6) & m_{11}=(-1,\,-1,\, 5) & m_{12}=(-1,\,-1,\, 4) \\
m_{13}=(-1,\,-1,\, 3) & m_{14}=(-1,\,-1,\, 2) & m_{15}=(-1,\,-1,\, 1) \\
m_{16}=(-1,\,-1,\, 0) & m_{17}=(-1,\, 0,\,-1) & m_{18}=(-1,\, 1,\,-1) \\
m_{19}=(0,\,-1,\,4) & m_{20}=(0,\, -1,\, -1)
\end{array}
\]
and $\{D_i\, | \, i\not=1,4,5 \}$ is linearly independent by system (\ref{linearRelation}). 
The lattice $L:=\langle r_* D_i\, | \, i\not=1,4,5 \rangle_\mathbb{Z}$ has $\rk L=17$, and $\discr L=2$ by an explicit calculation of its intersection matrix using formulas (\ref{SelfIntersection}) and (\ref{Intersection}). 
Since the discriminant of $L$ is square-free, $L$ is primitively embedded into the $K3$ lattice, and thus $L=\Pic(\Delta)$. 

The dual polytope $\Delta^*$ associates a toric $3$-fold with $6$ toric divisors $D'_i$ corresponding to  
\[
\begin{array}{lll}
v_1=(0,\,1,\,0) & v_2=(0,\,0,\,1) & v_3=(-3,\,-2,\,0) \\
v_4=(-5,\,-3,\,-1) & v_5=(1,\,0,\,0) & v_6=(-1,\,-1,\,0) 
\end{array}
\]
and $\{D'_i\, |\, i\not=1,2,3\}$ is linearly independent by system (\ref{linearRelation}). 
By re-taking the generators as $\{r_*\, D'_4,\, r_*\, D'_5-8r_*\, D'_4,\, r_*\, D'_6-r_*\, D'_4\}$, the lattice $L':=\langle r_*D'_i \, | \, i\not= 1,2,3 \rangle_\mathbb{Z}$ has an intersection matrix $U\oplus A_1$. 
Since $\discr L' = 2$ is square-free, $L'$ is primitively embedded into the $K3$ lattice, and thus $\Pic(\Delta^*)=L'=\left(\mathbb{Z}^3,\, U\oplus A_1\right)$. 

Similarly, $U\oplus\Pic(\Delta^*)$ is primitively embedded into the $K3$ lattice. 
Besides, $\discr (U\oplus \Pic(\Delta^*)) = -2=-\discr \Pic(\Delta)$, and moreover, $A_{\Pic(\Delta)}\simeq A_{U\oplus \Pic(\Delta^*)}\simeq \mathbb{Z}\slash 2\mathbb{Z}$. \\

\noindent
${\bf 6. }\, Q_{18}$ and $E_{30}$ \, 
The polytope $\Delta$ is given in subsection 4.13 of \cite{MU} and the associated toric $3$-fold has $19$ toric divisors $D_i$ corresponding to 
\[
\begin{array}{lll}
m_1=(1,\,-1,\,-1) & m_2=(-1,\,-1,\,-1) & m_3=(-1,\,-1,\,8) \\
m_4=(1,\,-1,\,0) & m_5=(-1,\,2,\,-1) & m_6=(0,\,-1,\,-1) \\
m_7=(-1,\,-1,\, 0) & m_8=(-1,\,-1,\, 1) & m_9=(-1,\,-1,\, 2) \\
m_{10}=(-1,\,-1,\, 3) & m_{11}=(-1,\,-1,\, 4) & m_{12}=(-1,\,-1,\, 5) \\
m_{13}=(-1,\,-1,\, 6) & m_{14}=(-1,\,-1,\, 7) & m_{15}=(0,\,-1,\,4) \\
m_{16}=(-1,\, 0,\,-1) & m_{17}=(-1,\, 1,\,-1) & m_{18}=(-1,\, 0,\, 5) \\
m_{19}=(-1,\, 0,\, 2)
\end{array}
\]
and $\{D_i\, | \, i\not=1,4,5 \}$ is linearly independent by system (\ref{linearRelation}). 
The lattice $L:=\langle  r_* D_i\, | \, i\not=1,4,5 \rangle_\mathbb{Z}$ has $\rk L=16$, and $\discr L=-3$ by an explicit calculation of its intersection matrix using formulas (\ref{SelfIntersection}) and (\ref{Intersection}). 
Since the discriminant of $L$ is square-free, $L$ is primitively embedded into the $K3$ lattice, and thus $L=\Pic(\Delta)$. 

The dual polytope $\Delta^*$ associates a toric $3$-fold with $7$ toric divisors $D'_i$ corresponding to 
\[
\begin{array}{llll}
v_1=(0,\,0,\,1) & v_2=(1,\,0,\,0) & v_3=(-4,\,-3,\,-1) & v_4=(-3,\,-2,\,0) \\
v_5=(0,\,1,\,0) & v_6=(-2,\,-1,\,0) & v_7=(-1,\,0,\,0) 
\end{array}
\]
and $\{D'_i\, |\, i\not=1,2,5\}$ is linearly independent by system (\ref{linearRelation}). 
By re-taking the generators as $\{r_*\, D'_3,\, r_*\, D'_4+r_*\, D'_3,\, r_*\, D'_6-r_*\, D'_3,\, r_*\, D'_7\}$, the lattice $L':=\langle r_*D'_i \, | \, i\not= 1,2,5 \rangle_\mathbb{Z}$ has an intersection matrix $U\oplus A_2$ . 
Since $\discr L'=-3$ is square-free, $L'$ is primitively embedded into the $K3$ lattice, and thus $\Pic(\Delta^*) = L'=\left(\mathbb{Z}^4,\, U\oplus A_2\right)$. 

Similarly, $U\oplus\Pic(\Delta^*)$ is primitively embedded into the $K3$ lattice. 
Besides, $\discr (U\oplus \Pic(\Delta^*)) = 3= -\discr \Pic(\Delta)$, and moreover, 
$A_{\Pic(\Delta)}\simeq A_{U\oplus \Pic(\Delta^*)}\simeq \mathbb{Z}\slash 3\mathbb{Z}$. 
\QED

\begin{remark} \rm
The choice of $\Delta$ is not actually unique. 
However, for any possible reflexive polytopes $\Delta$ for transpose-dual pairs outside Table \ref{NicePairs}, we have a relation $\rho(\Delta)+\rho(\Delta^*)\not= 20$. 
\end{remark}

\begin{thm} \label{MainResult}
The families $\mathcal{F}_\Delta$ and $\mathcal{F}_{\Delta^*}$ are lattice mirror symmetric for polytopes $\Delta$ in Corollary \ref{SurjectiveRestriction}, that is, $\Pic(\Delta)\simeq (U\oplus\Pic(\Delta^*))^\perp_{\Lambda_{K3}}$. 
Explicitely, the lattices $\Pic(\Delta)$ are given as in the Table \ref{NicePairs2}. 
\[
\begin{array}{cccccc} 
\textnormal{Singularity} & \Pic(\Delta) & \rho(\Delta) & \rho(\Delta^*) & \Pic(\Delta^*) & \textnormal{Singularity} \\ 
\hline
\hline
Q_{12} & U\oplus E_6\oplus E_8 & 16 & 4 & U\oplus A_2 & E_{18}\\
\hline
Z_{1,0} & U\oplus E_7\oplus E_8 & 17 & 3 & U\oplus A_1 & E_{19}\\
\hline
E_{20} & U\oplus E_8^{\oplus 2} & 18 & 2 & U & E_{20}  \\
\hline
Q_{2,0} & U\oplus A_6\oplus E_8 & 16 & 4 & U\oplus C_8^6  & Z_{17}\\
\hline
E_{25} & U\oplus E_7\oplus E_8 & 17 & 3 & U\oplus A_1 & Z_{19}\\
\hline
Q_{18} & U\oplus E_6\oplus E_8 & 16 & 4 & U\oplus A_2 & E_{30}\\
\hline
\end{array}
\]
\begingroup
\captionof{table}{Picard lattices for lattice mirror symmetric pairs}\label{NicePairs2}
\endgroup
\end{thm}
{\sc Proof.} \, 
For a lattice $L$, denote by $A_L$ the discriminant group of $L$, and $q_L$ the quadratic form of $L$. 
We see in Proposition \ref{ExplicitLattices} that 
\[
A_{\Pic(\Delta)} \simeq A_{U\oplus \Pic(\Delta^*)}, 
\textnormal{ and }
q_{\Pic(\Delta)} = -q_{U\oplus \Pic(\Delta^*)}
\]
for each case. 
Thus by Proposition 1.6.1 in \cite{Nikulin80}, there is an isometry 
\[
\Pic(\Delta) \simeq (U\oplus \Pic(\Delta^*))^\perp_{\Lambda_{K3}}. 
\]
So $\mathcal{F}_\Delta$ and $\mathcal{F}_{\Delta^*}$ are lattice mirror symmetric. 
We shall determine $\Pic(\Delta)$. \\

\noindent
${\bf 1. }\, Q_{12}$ and $E_{18}$ \, 
By Corollary \ref{NLem4.3}
\[
\Pic(\Delta) \simeq (U\oplus U\oplus A_2)^\perp_{\Lambda_{K3}} \simeq  
U\oplus (A_2)^\perp_{E_8}\oplus E_8 \simeq 
U\oplus E_6\oplus E_8. 
\]

\noindent
${\bf 2. }\, Z_{1,0}$ and $E_{19}$ \, 
By Corollary \ref{NLem4.3}
\[
\Pic(\Delta) \simeq (U\oplus U\oplus A_1)^\perp_{\Lambda_{K3}}  \simeq 
U\oplus (A_1)^\perp_{E_8} \oplus E_8 \simeq 
U\oplus E_7\oplus E_8. 
\]

\noindent
${\bf 3. }\, E_{20}$ and $E_{20}$ \, 
Since $\Lambda_{K3}\simeq U^{\oplus 3}\oplus E_8^{\oplus 2}$
\[
\Pic(\Delta)\simeq (U\oplus U)^\perp_{\Lambda_{K3}} \simeq U\oplus E_8^{\oplus 2}. 
\]

\noindent
${\bf 4. }\, Q_{2,0}$ and $Z_{17}$ \, 
By Corollary \ref{NLem4.3}
\[
\Pic(\Delta)\simeq (U\oplus U\oplus C_8^6)^\perp_{\Lambda_{K3}} \simeq 
U\oplus (C_8^6)^\perp_{E_8} \oplus E_8 \simeq 
U\oplus A_6\oplus E_8. 
\]

\noindent
${\bf 5. }\, E_{25}$ and $Z_{19}$ \, 
By Corollary \ref{NLem4.3}
\[
\Pic(\Delta)\simeq (U\oplus U\oplus A_1)^\perp_{\Lambda_{K3}} \simeq 
U\oplus (A_1)^\perp_{E_8} \oplus E_8 \simeq 
U\oplus E_7\oplus E_8. 
\]

\noindent
${\bf 6. }\, Q_{18}$ and $E_{30}$ \, 
By Corollary \ref{NLem4.3}
\[
\Pic(\Delta)\simeq (U\oplus U\oplus A_2)^\perp_{\Lambda_{K3}} \simeq  
U\oplus (A_2)^\perp_{E_8} \oplus E_8 \simeq 
U\oplus E_6\oplus E_8. 
\]
Thus the assertions are verified.
\QED

\section{Application} \label{Application}
Denote by $\Delta_B$ the reflexive polytope obtained in \cite{MU} for a singularity of type $B$. 
As is seen in Table \ref{NicePairs}, there are isometric Picard lattices  $\Pic(\Delta_{Q_{12}})\simeq \Pic(\Delta_{Q_{18}})$, and $\Pic(\Delta_{Z_{1,0}})\simeq \Pic(\Delta_{E_{25}})$. 
Families $\mathcal{F}_{\Delta_B}$ and $\mathcal{F}_{\Delta_D}$ are said to be {\it essentially the same} if general members in these families are birationally equivalent. 
Not only the Picard lattices are isometric, but also we have 
\begin{prop} \label{BirFamilies}
The families $\mathcal{F}_{\Delta_{Q_{12}}}$ $($resp. $\mathcal{F}_{\Delta_{Z_{1,0}}})$ and $\mathcal{F}_{\Delta_{Q_{18}}}$ $($resp. $\mathcal{F}_{\Delta_{E_{25}}})$ are essentially the same. 
\end{prop}
\textsc{Proof. }\, 
It is directly shown that the polytopes $\Delta_{Q_{12}}$ (resp. $\Delta_{Z_{1,0}}$) and $\Delta_{Q_{18}}$ (resp. $\Delta_{E_{25}}$) are isometric. 
Indeed, define an invertible matrix in $GL_3(\mathbb{Z})$ as 
\[
M_1 = \left( \begin{array}{ccc} -1 & 0 & 1 \\ 1 & 1 & -2 \\ 2 & -3 & 2 \end{array}\right)
\left( \textnormal{resp. } 
M_2 = \left( \begin{array}{ccc} -1 & 2 & 1 \\ 1 & -3 & -1 \\ 2 & 0 & -1 \end{array}\right) \right)
\]
and one obtains an isometry
\[
m \, M_1 = m' (\textnormal{resp. } m\, M_2 = m' )
\]
that sends $m\in\Delta_{Q_{18}}$ to $m'\in\Delta_{Q_{12}}$  (resp. $m\in\Delta_{E_{25}}$ to $m'\in\Delta_{Z_{1,0}}$). 
Therefore, there exists an explicit projective transformation that maps each Laurent polynomial in $|{-}K_{\mathbb{P}(\Delta_{Q_{12}})}|$ (resp. $|{-}K_{\mathbb{P}(\Delta_{Z_{1,0}})}|$) to a Laurent polynomial in $|{-}K_{\mathbb{P}(\Delta_{Q_{18}})}|$ (resp. $|{-}K_{\mathbb{P}(\Delta_{E_{25}})}|$). 
This mapping also birationally sends a general member in $\mathcal{F}_{\Delta_{Q_{12}}}$ (resp. $\mathcal{F}_{\Delta_{Z_{1,0}}}$) to a general member in  $\mathcal{F}_{\Delta_{Q_{18}}}$ (resp. $\mathcal{F}_{\Delta_{E_{25}}}$). 
Thus the statement is proved. 
\QED \\

We conclude our study to remark that general members in compactifications of non-equivalent singularities can be transformed via a reflexive polytope.

\hfill \textsc{Makiko Mase} \\ 
\hfill e-mail: mtmase@arion.ocn.ne.jp \\
\hfill {\small{Department of Mathematics and Information Sciences, Tokyo Metropolitan University}} \\
\hfill {\small{192-0397 1-1 Minami Osawa, Hachioji-shi, Tokyo, Japan. }}\\
\hfill {\small{Osaka City University Advanced Mathematical Institute}} \\
\hfill {\small{558-8585 3-3-138 Sugimoto-cho, Sumiyoshi-ku, Osaka, Japan. }}\\

\begin{thebibliography}{10}
\bibitem{Arnold75}
Arnol'd, V. I., 
Critical points of smooth functions and their normal forms, 
Russian Math. Surveys 30, 1--75 (1975). 
\bibitem{BHPV}
Barth, W.P. and Hulek, K. and Peters, C.A.M. and Van de Ven, A., 
{\it Compact Complex Surfaces, Second ed}, 
Springer (2004). 
\bibitem{DolgachevMirror}
Dolgachev, I., 
Mirror symmetry for lattice polarized {$K3$} surfaces,  
J. Math. Sci. 81, 2599--2630 (1996). 
\bibitem{EbelingPloog}
Ebeling, W. and Ploog, D., 
A geometric construction of {C}oxeter-{D}ynkin diagrams of bimodal singularities, 
Manuscripta Math. 140, 195--212 (2013). 
\bibitem{Iano-Fletcher}
Iano-Fletcher, A. R., 
Working with weighted complete intersections, 
in Explicit Birational Geometry of 3-folds, Alessio Corti and Miles Reid (eds.) London Mathematical Society Lecture Note Series 281, 101--173 (2000). 
\bibitem{Kobayashi}
Kobayashi, M., 
Duality of weights, mirror symmetry and {A}rnold's strange duality,
Tokyo J. Math., 31, 225--251 (2008). 
\bibitem{MU}
Mase, M.and Ueda, K, 
A note on bimodal singularities and mirror symmetry,
Manuscripta Math.(online 31 August 2014) 146, 153--177 (2015). 
\bibitem{Nikulin80} 
Nikulin, V.V., 
Integral symmetric bilinear forms and some of their applications, 
Math. USSR-Izv. 14, 103--167 (1980). 
\bibitem{Nishiyama96}
Nishiyama, K., 
The {J}acobian fibrations on some $K3$ surfaces and their {M}ordell-{W}eil groups, 
Japan J. Math. 22, 293--347 (1996). 
\bibitem{Oda78}
Oda, T., 
{\it Torus Embeddings and Applications}, 
Springer, Tata Institute of Fundamental Research Lectures 57 (1978). 
\bibitem{C3f}
Reid, M.,
Canonical $3$-folds, 
in Journe\'es de g\'eom\'etrie alg\'ebrique d'Angers, edited by A. Beauville, Sijthoff and Noordhoff, Alphen, 273--310 (1980). 
\bibitem{Yonemura}
Yonemura, T., 
Hypersurface simple $K3$ singularities, 
T\^ohoku Math. J. 42, 351--380 (1990). 
\end{thebibliography}
\end{document}